\documentclass[reqno,11pt]{amsart}
\usepackage{amsmath}
\usepackage{amssymb}
\usepackage{amsthm,amsxtra}
\usepackage{cite}
\usepackage{enumerate}
\usepackage[mathscr]{eucal}
\usepackage{float}
\usepackage[left=30mm,right=30mm,top=30mm,bottom=28mm]{geometry}
\usepackage{tikz-cd}
\usepackage{tkz-graph}
\usepackage{tkz-berge}
\usetikzlibrary{positioning,arrows,shapes.geometric,trees,matrix}
\usepackage{graphicx}
\usepackage{mathtools}
\usepackage{lineno}
\usetikzlibrary{positioning}
\usepackage{adjustbox}

\newtheorem{Th}{Theorem}[section]

\newtheorem{Lemma}{Lemma}[section]
\newtheorem{Prop}{Proposition}[section]

{}

\newcommand{\nb}{\mathbb{N}}

\newcommand{\rc}{\mathcal{R}}

\newcommand{\jc}{\mathcal{J}}
\newcommand{\ic}{\mathcal{I}}
\newcommand{\ac}{\mathcal{A}}
\newcommand{\bc}{\mathcal{B}}
\newcommand{\bk}{\mathfrak{B}}
\newcommand{\bb}{\mathfrak{b}}

\newcommand{\fc}{\mathcal{F}}

\newcommand{\oc}{\mathcal{O}}
\newcommand{\pc}{\mathcal{P}}
\newcommand{\pf}{\mathcal{P}}

\newcommand{\cs}{\mathcal{S}}

\newcommand{\uc}{\mathcal{U}}

\newcommand{\vc}{\mathcal{V}}
\newcommand{\wc}{\mathcal{W}}

\newcommand{\zc}{\mathcal{Z}}

\newcommand{\obs}{\oc_{\bk^s}}
\newcommand{\bs}{{\bk^s}}

\newcommand{\ofbs}{\oc_{f(\bk)^s}}
\newcommand{\gbs}{\gamma_{\bk^s}}

\newcommand{\Gbs}{\Gamma_{\bk^s}}
\newcommand{\IGbs}{\ic\text{-}\Gamma_{\bk^s}}
\newcommand{\IOGbs}{\ic_1\text{-}\Gamma_{\bk^s}}
\newcommand{\ITGbs}{\ic_2\text{-}\Gamma_{\bk^s}}
\newcommand{\igbs}{\ic\text{-}\gamma_{\bk^s}}
\newcommand{\iogbs}{\ic_1\text{-}\gamma_{\bk^s}}
\newcommand{\itgbs}{\ic_2\text{-}\gamma_{\bk^s}}
\newcommand{\JGbs}{\jc\text{-}\Gamma_{\bk^s}}
\newcommand{\JOGbs}{\jc_1\text{-}\Gamma_{\bk^s}}
\newcommand{\JTGbs}{\jc_2\text{-}\Gamma_{\bk^s}}
\newcommand{\jgbs}{\jc\text{-}\gamma_{\bk^s}}
\newcommand{\jogbs}{\jc_1\text{-}\gamma_{\bk^s}}
\newcommand{\jtgbs}{\jc_2\text{-}\gamma_{\bk^s}}
\newcommand{\ibs}{\ic\text{-}\bk^s}
\newcommand{\iobs}{\ic_1\text{-}\bk^s}
\newcommand{\itbs}{\ic_2\text{-}\bk^s}
\newcommand{\IGfbs}{\ic\text{-}\Gamma_{f(\bk)^s}}
\newcommand{\JGfbs}{\jc\text{-}\Gamma_{f(\bk)^s}}
\newcommand{\igfbs}{\ic\text{-}\gamma_{f(\bk)^s}}
\newcommand{\jgfbs}{\jc\text{-}\gamma_{f(\bk)^s}}
\newcommand{\K}{K}
\newcommand{\KB}{KB}
\DeclareMathOperator{\Fin}{Fin}
\DeclareMathOperator{\Sf}{S_{fin}}
\DeclareMathOperator{\S1}{S_1}

\DeclareMathOperator{\Gf}{G_{fin}}
\DeclareMathOperator{\G1}{G_1}
\DeclareMathOperator{\sle}{sl_e}
\DeclareMathOperator{\slt}{sl_t}
\DeclareMathOperator{\fupc}{FUPC}
\DeclareMathOperator{\fup}{FUP}
\DeclareMathOperator{\Split}{Split}
\newcommand{\bp}{\begin{proof}}
\newcommand{\ep}{\end{proof}}

\begin{document}

\title[Certain observations on bornological covers using ideals]{Certain observations on selection principles related to bornological covers using ideals}
\author[D. Chandra, P. Das  and S. Das]{Debraj Chandra$^{\dag}$, Pratulananda Das$^*$ and Subhankar Das$^*$\ }
\address{\llap{$\dag$\,}Department of Mathematics, University of Gour Banga, Malda-732103, West Bengal, India}
\email{debrajchandra1986@gmail.com}

\address{\llap{*\,}Department of Mathematics, Jadavpur University, Kolkata-700032, West Bengal, India}
\email{pratulananda@yahoo.co.in, subhankarjumh70@gmail.com}

\thanks{The second author (as PI) and the third author (as Fellow) are thankful to NBHM for the research project (No. 02011/9/2022 / NBHM (RP)/R\&D II/ 10387) during the tenure of which this work has been done.}
\subjclass{Primary: 54D20; Secondary: 54C35, 54A25 }

\begin{abstract}
We study selection principles related to bornological covers using the notion of ideals. We consider ideals $\ic$ and $\jc$ on $\omega$ and standard ideal orderings $\KB, \K$. Relations between cardinality of a base of a bornology with certain selection principles related to bornological covers are established using cardinal invariants such as modified pseudointersection number, the unbounding number and slaloms numbers. When $\ic\leq_\square \jc$ for ideals $\ic, \jc$ and $\square\in \{1\text{-}1,\KB,\K\}$, implications among various selection principles related to bornological covers are established. Under the assumption that ideal $\ic$ has a pseudounion we show equivalences among certain selection principles related to bornological covers. Finally, the $\ibs$-Hurewicz property of $X$ is investigated. We prove that $\ibs$-Hurewicz property of $X$ coincides with the $\bs$-Hurewicz property of $X$ if $\ic$ has a pseudounion. Implications or equivalences among selection principles, games and $\ibs$-Hurewicz property which are obtained from our investigations are described in diagrams.
\end{abstract}
\maketitle
\smallskip
\noindent{\bf\keywordsname{}:} {Bornology, open $\bs$-cover, $\gamma_\bs$-cover, $\igbs$-cover, selection principles, topology of strong uniform convergence, ideal, pseudounion, ideal ordering.}


\section{Introduction}
We follow the notations and terminologies of \cite{arh, engelking, hh, mcnt}. The Selection Principle Theory is a widely known discipline of Mathematics having a long history dated back to 1920-30's. Scheepers in his celebrated papers \cite{cooc1,cooc2} had systematically formulated selection principles related to open covers. Since then this subject has been developed by many well known researchers. One can consult the survey papers \cite{sur1, sur2,sur3} to get a full view of the study of the Selection Principle Theory. In this paper, we study selection principles related to open covers from the bornological perspective using the notion of ideals. We first recall that a \textit{bornology} \cite{hh} $\bk$ on a space $X$ is a collection of subsets of $X$ which satisfies the following conditions: $(1)$ it is hereditary, i.e., for $B\in \bk$ and $B^\prime\subseteq B$, $B^\prime\in \bk$, $(2)$ it is closed under taking finite unions, i.e., for $B,B^\prime\in \bk$, $B\cup B^\prime\in \bk$, $(3)$ it is a cover of $X$, i.e., for $x\in X$ there exists a $B\in \bk$ with $x\in B$. A \textit{base} $\bk_0$ for $\bk$ is a subcollection of $\bk$ such that for $B\in \bk$ there exists a $B_0\in \bk_0$ such that $B\subseteq B_0$. The notion of a base for a bornology is  very useful for our bornological investigations of selection principles related to open covers. There is a connection between a base of a bornology with certain selection principles related to bornological covers. When members of a base are closed (compact) it is called closed (compact) base. The family $\fc$ of all finite subsets of $X$ and the family $\pc(X)$ of all subsets of $X$ respectively are bornologies on $X$ (for more examples see \cite{dcpdsd}).

The notion of bornology plays an important role in Functional Analysis, Topology, Selection Principles Theory etc. As far as bornological investigations to Selection Principles Theory is concerned, it has been started when the notions of \textit{strongly uniformly continuous function} and the \textit{topology of strong uniform convergence on a bornology} have been introduced by Beer and Levi \cite{bl} in 2009. Since then the study of selection principles related to open covers from bornological perspective has been started. One can see papers \cite{bl, cmh, cmk, dcpdsd, dcpdsd2, pddcsd3} for more informations. In \cite{sddc} the authors have taken a general approach to study selection principles related to bornological covers using the notion of statistical convergence in metric spaces. Generalized notions of certain bornological covers ($\igbs$-cover) and some bornological properties ($\ibs$-Hurewicz property etc.) have been investigated. In \cite{sd} the notion of ideals has been used to study selection principles related to bornological covers..

A family $\ic$ of subsets of $\pc(\omega)$ is called an \textit{ideal} if the following conditions hold: $(1)$ it is a hereditary family, i.e., $A\in \ic$ and $B\subseteq A$ implies $B\in \ic$, $(2)$ it is closed under finite unions, i.e., $A\cup B\in \ic$ for any $A,B\in \ic$, $(3)$ it contains all finite subsets of $\omega$ and $\omega\not\in \ic$. $\ic$ is also called a \textit{proper} and \textit{admissible} ideal. The set of all finite subsets of $\omega$ is the \textit{smallest} ideal and it is denoted by $\Fin$. For a family $\ac$ subset of $\pc(\omega)$ $\ac^\ast=\{\omega\setminus A:A\in \ac\}$. $\ac$ is called a \textit{filter} if $\ac^\ast$ is an ideal. Also $\ic^+=\pc(\omega)\setminus \ic$. For recent developments on the ideal approach to selection principles theory and related covers, see \cite{vsjs, js(QN), js, lbpdjs,sbjslz,amjs} .

The main objective of the paper is to study selection principles related to bornological covers using the notions of ideals. Unlike the paper \cite{sd} where a single general ideal $\ic$ on $\omega$ has been used, in this paper more than one ideals ($\ic, \jc$) are considered and standard orderings ($\KB,\K$) on ideals are applied for the investigations. This reveals more information about behaviour of selection principles related to bornological covers. More importantly we use the concept of pseudounion on $\ic$.  When an ideal $\ic$ has a pseudounion we show that every $\igbs$-cover of $X$ contains a subset which is a $\gbs$-cover of $X$. These results play a crucial role to show equivalences among various selection principles related to bornological covers and certain bornological properties which were otherwise previously unknown.

The paper is structured as follows. In Section $3$, we begin with some observations on interactions between cardinality of a base of a bornology with various selection principles related to bornological covers. Cardinal invariants such as modified pseudointersection number, the unbounding number, slaloms etc. are considered. Behaviours of certain selection principles related to bornological covers under continuous images are presented. Considering ideals $\ic$ and $\jc$ on $\omega$ and applying the standard ideal ordering $\KB$ and $\K$, implications among various selection principles related to bornological covers are established and these are described in diagrams (Figures~\ref{diag1},\ref{diag3},\ref{diag2}). Assuming that $\ic$ has a pseudounion, we show equivalences among certain selection principles related to bornological covers (Figure~\ref{diag4}). Some observations on splittability are presented (Figure~\ref{diag7}). Section $4$ deals with the $\ibs$-Hurewicz property of $X$. We proved that the $\iobs$-Hurewicz property of $X$ implies the $\itbs$-Hurewicz property of $X$, where $\ic_1\leq_\square \ic_2$ with $\square\in \{1\text{-}1,\KB,\K\}$. Also the $\ibs$-Hurewicz property of $X$ coincides with the $\bs$-Hurewicz property if $\ic$ has a pseudounion. Finally, relations among the $\ibs$-Hurewicz property, $\Sf(\obs,\ic\text{-}\obs^{gp})$ and certain games related to bornological covers are presented (Figure~\ref{diag6}).

\section{Preliminaries}
Let $(X,d)$ be an infinite metric space and let $\omega$ denote the set of all natural  numbers.
For two nonempty classes of sets $\ac$ and $\bc$ of an infinite set $S$, we recall that

\noindent $\S1(\ac,\bc)$: For every sequence $\{A_n:n\in \omega \}$ of elements of $\ac$, there is a sequence $\{a_n:n\in \omega\}$ with $a_n\in A_n$ for each $n$ such that $\{a_n:n\in \omega\}\in \bc$ (\cite{cooc1,cooc2}).

\noindent $\Sf(\ac,\bc)$: For every sequence $\{A_n:n\in \nb\}$ of elements of $\ac$, there is a sequence $\{B_n:n\in \nb\}$ of finite (possibly empty) sets with $B_n\subseteq A_n$ for each $n$ such that $\bigcup_{n\in \nb}B_n\in \mathcal B$ (\cite{cooc1,cooc2}).

Corresponding to these selection principles there are infinitely long games.

\noindent $\G1(\ac,\bc)$ denotes the game played by two players ONE and TWO who play a round for each positive integer $n$. In the $n$-th round ONE chooses a set $A_n$ from $\ac$ and then TWO chooses an element $a_n\in  A_n$. TWO wins the play $\{A_1,a_1, \dotsc, A_n, a_n, \dotsc \}$ if $\{a_n :n\in \nb\}\in \bc$. Otherwise ONE wins.

\noindent $\Gf(\ac,\bc)$ denotes the game played by two players ONE and TWO who play a round for each positive integer $n$. In the $n$-th round ONE chooses a set $A_n$ from $\ac$ and then TWO chooses a finite (possibly empty) set $B_n\subseteq A_n$. TWO wins the play $\{A_1,B_1, \dotsc,A_n,B_n, \dotsc  \}$ if $\bigcup_{n\in \nb}B_n\in \bc$. Otherwise ONE wins.\\

For two sets $A, B$ $A$ is said to be \textit{almost contained} in $B$ denoted by $A\subseteq^\ast B$ if $A\setminus B$ is finite. Let $\ac$ be a family of subsets of $\omega$. A subset $B$ of $\omega$ is a \textit{pseudounion} of $\ac$ if $X\setminus B$ is infinite and $A\subseteq^\ast B$ for every $A\in \ac$ \cite{lbpdjs}. $\ac$ is called \textit{tall} if $\ac$ has no pseudounion. $\ac$ is said to have the finite union property (in short, $\fup$) if for any $A_1,\dotsc A_k\in \ac$, the set $\omega\setminus \cup_{i=1}^k A_i$ is infinite \cite{lbpdjs}. One can easily verify that a family $\ac$ has the $\fup$ if there is an ideal $\ic$ such that $\ac\subseteq \ic$. We say that $\ac$ has $(\ast)$ if $\ac$ has the finite union property, and $\ac$ has $(\ic)$ if $\ac\subseteq \ic$. Throughout $^{X}Y$ will stand for the set of all functions from $X$ to $Y$.

Let $\ic$ be an ideal on $\omega$. For $f,g\in {^\omega}\omega$ by $f\leq_\ic g$ we mean $\{n\in\omega:g(n)<f(n)\}\in \ic$ and a subset $A$ of ${^\omega}\omega$ is \textit{$\ic$-bounded} if there exists a function $g$ satisfying $f\leq_\ic g$ for all $f\in A$. The symbol $\bb_\ic$ denotes the unbounding number of $<{^\omega}\omega,\leq_\ic>$ and is defined as the minimal size of a $\leq_\ic$-unbounded family \cite{bfls}.

Let $\square\in \{1\text{-}1, \KB, \K \}$. A function $\phi:\omega\rightarrow \omega$ is called a $\square$-function if $\phi$ is a one-one, a finite-to-one and a arbitrary function respectively. For a function $\phi: \omega\rightarrow \omega$, $\ac_1\leqslant_\phi\ac_2$ implies that $\phi^{-1}(A)\in \ac_2$ for any $A\in \ac_1$. $\ac_1\leqslant_\square\ac_2$ implies that there exists a $\square$-function $\phi$ with $\ac_1\leqslant_\phi\ac_2$. Now for $\square\in \{1\text{-}1, \KB, \K \}$, an ideal $\jc$ of $\omega$ and $\triangle\in \{\ic,\ast\}$,
$\pf_\square(\triangle,\jc)=\min\{|\ac|:\ac\subseteq\pc(\omega) \text{ has } (\triangle)\wedge\ac\nleq_\square \jc\}$. Another formulation of $\pf_\square(\triangle,\jc)$ is
\[ \pf_\square(\triangle,\jc)=\min\{|\ac|:\ac\subseteq \pc(\omega) \text{ has } (\triangle)\hspace{.05in}\wedge\hspace{.05in} (\forall \phi\in {^\omega}\omega,\square\hspace{.05in} \exists\hspace{.05in} A\in \ac \text{ such that } \{n:\phi(n)\in A\}\in \jc^+)\}.\]
$\pf_\square(\jc)=\min\{\pf_\square(\ic,\jc):\ic \text{ is an ideal }\}$.\\

$[\pc,\rc]_\square$: For every sequence $\{V_n:n\in \omega\}\in \pc$ there exists a $\square$-function $\phi: \omega\rightarrow \omega$ such that $\{V_{\phi(n)}:n\in \omega\}\in \rc$, where $\square\in \{1\text{-}1, \KB, \K \}$.\\

For ideals $\ic,\jc\subseteq \pc(\omega)$ one can consider the following cardinal number \cite{js(QN)}.
\[ \lambda(\ic,\jc)=\min\{|\ac|:\ac\subseteq {^\omega}\ic\hspace{.05in}\wedge\hspace{.05in} (\forall \phi\in {^\omega}\omega\hspace{.05in} \exists\hspace{.05in} \{B_n:n\in \omega\}\in \ac \text{ such that } \{n:\phi(n)\in B_n\}\in \jc^+)\} \]
Let $\cs\subset {^\omega}\pc(\omega)$. Elements of $\cs$ are called \textit{slaloms}. $\cs$ has finite union property coordinatewisely (in short, $\fupc$) if the family $\{s(n):s\in \cs\}$ has the finite union property for each $n\in \omega$. By $\phi^{-1}(s)$ we mean $\{n:\phi(n)\in s(n)\}$ for $\phi\in {^\omega}\omega$ and $s\in \cs$. By $\phi\in^\jc$ we mean $\phi^{-1}(s)\in \jc^\ast$. We say that $\cs$ has $(\ast)$ if $\cs$ has the finite union property coordinatewisely. $\cs$ has $(\ic)$ if $\cs\subseteq ^{\omega}\ic$. We now state definitions of slaloms cardinal invariants (see \cite{js} for more details).
\[ \sle(\ic,\jc)=\min\{|\cs|:\cs\subseteq {^\omega}\ic\hspace{.05in} \wedge \hspace{.05in} (\forall \phi\in {^\omega}\omega\hspace{.05in} \exists\hspace{.05in} s\in \cs \text{ such that } \{n:\phi(n)\in s(n)\}\in \jc^+)\} \]
\[ \slt(\ic,\jc)=\min\{|\cs|:\cs\subseteq {^\omega}\ic \hspace{.05in}\wedge \hspace{.05in}(\forall  \phi\in {^\omega}\omega\hspace{.05in} \exists\hspace{.05in} s\in \cs \text{ such that } \phi\in^\jc s)\}\]
\[ \sle(\ast, \jc)=\min\{|\cs|:\cs\subseteq {^\omega}\pc(\omega) \text{ has }\fupc \hspace{.05in}\wedge \hspace{.05in} (\forall  \phi\in {^\omega}\omega\hspace{.05in} \exists\hspace{.05in} s\in \cs \text{ such that }\{n:\phi(n)\in s(n)\}\in \jc^+)\} \]
\[ \sle(\ast, \jc)=\min\{|\cs|:\cs\subseteq {^\omega}\pc(\omega) \text{ has } \fupc \hspace{.05in}\wedge \hspace{.05in}(\forall \phi\in {^\omega}\omega\hspace{.05in} \exists\hspace{.05in} s\in \cs \text{ such that } \phi\in^\jc s)\}\]

Let $\bk$ be a bornology on metric space $(X,d)$ with closed base. Let $B^\delta= \bigcup_{x\in B}S(x,\delta)$ for $B\in \bk$, $\delta>0$ and $S(x,\delta)=\{y\in X:\rho(x,y)<\delta\}$. For $B\in \bk$ and $\delta>0$ it is easy to see that $\overline{B^\delta}\subseteq B^{2\delta}$. We now recall definitions of some classes of bornological covers of $X$. Let $\uc$ be a cover of $X$ and $X\not\in \uc$. $\uc$ is a \textit{strong-$\bk$-cover} (in short, $\bs$-cover)\cite{cmh} if for $B\in \bk$ there are $U\in \uc$ and $\delta>0$ satisfying $B^\delta\subseteq U$. We say $\uc$ an \textit{open $\bs$-cover} if the members of $\uc$ are open. By $\oc_\bs$ we denote the collection of all open $\bs$-covers of $X$. $X$ is $\bs$-Lindel\"of \cite{cmk} if every $\bs$-cover of $X$ contains a countable subset which is a $\bs$-cover of $X$. An open cover $\{U_n:n\in \omega\}$ is a $\gbs$-cover \cite{cmh} (see also \cite{cmk}) of $X$ if it is infinite and for every $B\in \bk$ there are a $n_0\in \omega$ and a sequence $\{\delta_n:n\geq n_0\}$ of positive real numbers satisfying $B^{\delta_n}\subseteq U_n$ for all $n\geq n_0$. By $\Gbs$ we denote the collection of all $\gbs$-covers of $X$. An open cover $\uc$ of $X$ is \textit{$\bk^s$-groupable} \cite{dcpdsd} if there exists a sequence $\{\uc_n:n\in \omega\}$ of finite pairwise disjoint sets such that $\uc=\cup_{n\in \omega}\uc_n$ and for each $B\in \bk$ there are a $n_0\in \omega$ and a sequence $\{\delta_n:n\ge  n_0\}$ of positive real numbers with $B^{\delta_n}\subseteq U$ for some $U\in \uc_n$ for all $n\ge  n_0$. By $\obs^{gp}$ we denote the collection of all open $\bs$-groupable covers of $X$. $X$ is said to have the \textit{$\bs$-Hurewicz property} \cite{dcpdsd} if for each sequence $\{\uc_n:n\in \nb\}$ of open ${\bk^s}$-covers of $X$ there is a sequence  $\{\vc_n:n\in \nb\}$, where $\vc_n\subseteq \uc_n$ a finite set for each $n\in \omega$ such that for every $B\in \bk$ there are a $n_0\in \omega$ and a sequence $\{\delta_n:n\ge  n_0\}$ of positive real numbers satisfying $B^{\delta_n}\subseteq U$ for some $U\in \vc_n$ for all $n\ge  n_0$. The \textit{strong $\bk$-Hurewicz game} (in short, $\bk^s$-Hurewicz game) on $X$ is an infinitely long game played by two players ONE and TWO. In the $n$-th inning ONE chooses an open $\bk^s$-cover $\uc_n$ of $X$, then TWO chooses a finite set $\vc_n\subseteq \uc_n$. TWO wins the play: $\uc_1, \vc_1, \uc_2, \vc_2, \dotsc ,\uc_n, \vc_n, \dotsc$ if for each $B\in \bk$ there are a $n_0\in \omega$ and a sequence $\{\delta_n:n\ge  n_0\}$ of positive real numbers satisfying $B^{\delta_n}\subseteq U$ for some $U\in \vc_n$ for all $n\ge  n_0$. Otherwise ONE wins.

Let $\ic$ be an ideal on $\omega$. We now recall ideal variants of some bornological notions. A countable open cover $\{U_n:n\in \nb\}$ of $X$ is an \textit{$\ic$-$\gbs$-cover} if for $B\in \bk$ there is a sequence $\{\delta_n:n\in \omega\}$ of positive real numbers satisfying $\{n\in \nb:B^{\delta_n}\nsubseteq U_n\}\in \ic$ \cite{sd}. $X$ is said to have the \textit{$\ic$-strong-$\bk$-Hurewicz property} (in short, $\ic$-$\bs$-Hurewicz property) if for each sequence $\{\uc_n:n\in \omega\}$ of open $\bs$-covers of $X$, there is a sequence $\{\vc_n:n\in \omega\}$, where $\vc_n\subseteq \uc_n$ a finite set for each $n\in \omega$ such that for every $B\in \bk$ there is a sequence $\{\delta_n:n\in \omega\}$ of positive real numbers satisfying $\{n\in \nb:B^{\delta_n}\nsubseteq U$ for any $U\in \vc_n\}\in \ic$ \cite{sd}. The \textit{$\ic$-strong-$\bk$-Hurewicz game} (in short, $\ic$-$\bk^s$-Hurewicz game) on $X$ is an infinitely long game played by two players ONE and TWO. In the $n$-th inning ONE chooses an open $\bk^s$-cover $\uc_n$ of $X$. Then TWO chooses a finite set $\vc_n\subseteq \uc_n$. TWO wins the play: $\uc_1$, $\vc_1$, $\uc_2$, $\vc_2, \dotsc$, $\uc_n$, $\vc_n$, $\dotsc$ if for each $B\in \bk$ there is a sequence $\{\delta_n:n\in \omega\}$ of positive real numbers satisfying $\{n\in \nb:B^{\delta_n}\nsubseteq U$ for any $U\in \vc_n\}\in \ic$. Otherwise ONE wins \cite{sd}. An open cover $\uc$ of $X$ is \textit{$\ic$-strong-$\bk$-groupable} (in short, $\ic$-$\bs$-groupable) if there exists a sequence $\{\uc_n:n\in \omega\}$ of finite pairwise disjoint sets such that $\uc=\cup_{n\in \omega}\uc_n$ and for each $B\in \bk$ there is a sequence $\{\delta_n:n\in \nb\}$ of positive real numbers satisfying $\{n\in \nb:B^{\delta_n}\nsubseteq U$ for any $U\in \uc_n\}\in \ic$ \cite{sd}. By $\ic$-$\obs^{gp}$ we denote the collection of all open $\ibs$-groupable covers of $X$.

\section{Some observations on Selection Principles}

\subsection{Basic observations on selection principles}
The first Lemma is related to $\igbs$-cover (see \cite[Lemma 2.1]{dcpdsd2}).
\begin{Lemma}
\label{LIf-1}
Let $\bk$ be a bornology on a metric space $X$ with a compact base $\bk_0$. Let $f:X\rightarrow Y$ be a continuous function. If $\{U_n:n\in \omega\}$ is an $\ic$-$\gamma_{f(\bk)^s}$-cover of $X$, then $\{f^{-1}(U_n):n\in \omega\}$ is an $\ic$-$\gbs$-cover of $X$.
\end{Lemma}
\bp
Let $B\in \bk_0$. For $f(B)\in f(\bk)$ there is a sequence $\{\varepsilon_n:n\in \omega\}$ of positive real numbers such that $\{n:f(B)^{\varepsilon_n}\nsubseteq U_n\}\in \ic$. As $B$ is compact, $f$ is strongly uniformly continuous on $B$. For $\varepsilon_n>0$ there is a $\delta_n>0$ such that $f(B^{\delta_n})\subseteq f(B)^{\varepsilon_n}$. It is easy to see that $\{n:B^{\delta_n}\nsubseteq f^{-1}(U_n)\}\subseteq \{n:f(B)^{\varepsilon_n}\nsubseteq U_n\}\in \ic$. Hence $\{f^{-1}(U_n):n\in \omega\}$ is an $\ic$-$\gbs$-cover of $X$.
\ep

\begin{Prop}
\label{PIf-1}
Let $\bk$ be a bornology on a metric space $X$ with a compact base $\bk_0$ and let $\ic,\jc$ be ideals on $\omega$. Let $f:X\rightarrow Y$ be a continuous function. The following statements hold.\\
\noindent$(1)$ If $X$ satisfies $\S1(\IGbs,\JGbs)$, then $f(X)$ satisfies $\S1(\IGfbs,\JGfbs)$.\\
\noindent$(2)$ If $X$ satisfies $\S1(\obs,\JGbs)$, then $f(X)$ satisfies $\S1(\ofbs,\JGfbs)$.
\end{Prop}
\bp
We prove only $(1)$. Let $\{\uc_n:n\in \omega\}$ be a sequence of $\igfbs$-covers of $X$ and $\uc_n=\{U_{n,m}:m\in \omega\}$ for each $n$. Let $\vc_n=\{f^{-1}(U_{n,m}):m\in \omega\}$. By Lemma \ref{LIf-1}, $\vc_n$ is an $\igbs$-cover of $X$. By the given condition, there is a $m_n$ for each $n$ such that $\{f^{-1}(U_{n,m_n}):m\in \omega\}$ is a $\jgbs$-cover. We show that $\{U_{n,m_n}:m\in \omega\}$ is a $\jgfbs$-cover. Let $B^\prime\in f(\bk)$ and $B^\prime=f(B)$ where $B\in \bk$. Now for the aforementioned $B\in \bk$, we can choose a sequence $\{\delta_n:n\in \omega\}$ of positive real numbers such that $\{n:B^{\delta_n}\nsubseteq f^{-1}(U_{n,m_n})\}(=S)\in \ic$. If $n\in S$, then $B^{\delta_n}\nsubseteq f^{-1}(U_{n,m_n})$ and so $f(B^{\delta_n})\nsubseteq U_{n,m_n}$. Choose a $\eta_n>0$  so that $f(B)^{\eta_n}\nsubseteq U_{n,m_n}$. Now $n\not\in S$ implies $B^{\delta_n}\subseteq f^{-1}(U_{n,m_n})$ and subsequently $f(B)\subseteq U_{n,m_n}$. As $f(B)$ is compact, there is a $\varepsilon_n>0$ such that $f(B)^{\varepsilon_n}\subseteq U_{n,m_n}$. Define $\sigma(n)=\varepsilon_n$ if $n\not\in S$ and $\sigma(n)=\eta_n$ if $n\in S$. Clearly $\{n:f(B)^{\sigma_n}\nsubseteq U_{n,m_n}\}=S\in \jc$. Hence $f(X)$ satisfies $\S1(\IGfbs,\JGfbs)$.
\ep

\begin{Prop}
\label{PIf-2}
Let $\bk$ be a bornology on a metric space $X$ with a compact base $\bk_0$ and let $\ic,\jc$ be ideals on $\omega$. If $X$ satisfies $\S1(\IGbs,\JGbs)$, then every continuous image of $X$ into the Baire space $\omega^\omega$ is $\jc$-bounded.
\end{Prop}
\bp
Let $\phi:X\rightarrow \omega^\omega$ be a continuous function. By Proposition \ref{PIf-1}(1), $\phi(X)$ satisfies $\S1(\ic\text{-}\Gamma_{\phi(\bk)^s},\jc\text{-}\Gamma_{\phi(\bk)^s})$. Let $U^n_k=\{f\in ^{\omega}\omega:f(n)\leq k\}$ and $\uc_n=\{U_{n_k}:k\in \omega\}$. For each $n$ one can show that $\uc_n$ is an $\igbs$-cover of $\phi(X)$. Applying $\S1(\ic\text{-}\Gamma_{\phi(\bk)^s},\jc\text{-}\text{-}\Gamma_{\phi(\bk)^s})$ to $\{\vc_n:n\in \omega\}$, we can choose a $k_n$ for each $n$ such that $\{U_{n,{k_n}}:n\in \omega\}$ is a $\jgbs$-cover. Define $h:\omega\rightarrow \omega$ such that $h(n)=k_n$. Let $f\in \phi(X)$. We claim that the set $S=\{n:h(n)<f(n)\}\in \jc$. Consider $B=\{f\}\in \phi(\bk)$. Then there is a sequence $\{\delta_n:n\in \omega\}$ of positive real numbers such that $\{n:B^{\delta_n}\nsubseteq U_{n,{k_n}}\}\in \jc$. Now take  $n\in S$ and observe that $h(n)<f(n)$ which means $k_n<f(n)$. Consequently $f\not\in U_{n,{k_n}}$ and $S$ being contained in $\{n:B^{\delta_n}\nsubseteq U_{n,{k_n}}\}$ is a member of $\jc$. Hence $\phi(X)$ is $\jc$-bounded.
\ep

Before moving forward we state some observations related to subset schema $\binom{\ac} {\bc}$. The relations among $\binom{\ac} {\bc}$, where $\ac,\bc\in \{\obs, \Gbs,\IGbs,\JGbs\}$ and $\ic,\jc$ be any two ideals, is described in Figure \ref{diag1}.

\begin{figure}[h]
 \begin{adjustbox}{keepaspectratio,center}
\begin{tikzcd}[column sep=2ex,row sep=5ex,arrows={crossing over}]
\binom{\IGbs} {\Gbs}\arrow[r]&\binom{\IGbs} {\JGbs}\\
\binom{\obs} {\Gbs}\arrow[r]\arrow[u]&\binom{\obs} {\JGbs}\arrow[u]
\end{tikzcd}
\end{adjustbox}
 \caption{Diagram of the selection principles $\binom{\ac} {\bc}$}
 \label{diag1}
\end{figure}

\begin{Th}
\label{TI-2}
Let $\bk$ be a bornology on a metric space $X$ with closed base and let $\ic$ be an ideal on $\omega$. Then the following statements hold.\\
\noindent$(1)$ $X$ satisfy $\binom{\obs} {\IGbs}$ and $\S1(\IGbs,\Gbs)$ if and only if $X$ satisfies $\S1(\obs,\Gbs)$.\\
\noindent$(2)$ $X$ satisfy $\binom{\obs} {\IGbs}$ and $\S1(\IGbs,\IGbs)$ if and only if $X$ satisfies $\S1(\obs,\IGbs)$.\\
\noindent$(3)$ $X$ satisfy $\binom{\IGbs} {\Gbs}$ and $\S1(\Gbs,\Gbs)$ if and only if $X$ satisfies $\S1(\IGbs,\Gbs)$.
\end{Th}
%

\begin{Th}
\label{TI-3}
Let $\bk$ be a bornology on a metric space $X$ with closed base and let $\ic, \jc$ be ideals on $\omega$. Then the following statements hold.\\
\noindent$(1)$ $X$ satisfy $\binom{\IGbs} {\Gbs}$ and $\S1(\obs,\IGbs)$ if and only if $X$ satisfies $\S1(\obs,\Gbs)$.\\
\noindent$(2)$ $X$ satisfy $\binom{\JGbs} {\Gbs}$ and $\S1(\IGbs,\JGbs)$ if and only if $X$ satisfies $\S1(\IGbs,\Gbs)$.
\end{Th}
%

\subsection{Results related to cardinality}
We now present observations under certain cardinality constraints of a base of a bornology.
\begin{Th}
\label{TI-1}
Let $\bk$ be a bornology on a metric space $X$ with a closed base $\bk_0$ and let $\ic, \jc$ be ideals on $\omega$. The following statements hold.\\
\noindent$(1)$ If $|\bk_0|<\pf_\square(\ic,\jc)$, then $[\IGbs,\JGbs]_\square$ holds.\\
\noindent$(2)$ If $|\bk_0|<\pf_\square(\jc)$, then $[\obs,\JGbs]_\square$ holds, provided $X$ is $\bs$-Lindel\"{o}f.
\end{Th}
\bp
$(1)$. Let $\{V_n:n\in \omega\}$ be an $\igbs$-cover of $X$. For a $B\in \bk_0$ there is a sequence $\{\delta_n:n\in \omega\}$ of positive real numbers satisfying $\{n\in \omega:B^{\delta_n}\nsubseteq V_n\}\in \ic$. Let $A_B=\{n\in \omega:B^{\delta_n}\nsubseteq V_n\}$. Consider the collection $\ac_{\bk_0}=\{A_B:B\in \bk_0\}$. Clearly $|\ac_{\bk_0}|\leq |\bk_0|$ and therefore $\ac_{\bk_0}\leq_\square \jc$. There exists a $\square$-function $\phi:\omega\rightarrow \omega$ such that $\phi^{-1}(A_B)\in \jc$ for all $A_B\in \ac_\bk$. We claim that $\{V_{\phi(m)}:m\in \omega\}$ is a $\jgbs$-cover of $X$. Let $B\in \bk_0$. We have  $\phi^{-1}(A_B)\in \jc$. Define $\sigma_m=\delta_{\phi(m)}$ for $m\in \omega$. The proof will be complete if we can show that $\{m:B^{\sigma_m}\nsubseteq V_{\phi(m)}\} (=C)\in \jc$.
If $m\in C$, then $B^{\sigma_m}\nsubseteq V_{\phi(m)}$ that is $B^{\delta_{\phi(m)}}\nsubseteq V_{\phi(m)}$ which implies $\phi(m)\in A_B$ and so $m\in \phi^{-1}(A_B)$. Hence $C\subseteq \phi^{-1}(A_B)\in \jc$.

$(2)$. Let $\{U_n:n\in \omega\}$ be an open $\bs$-cover of $X$. For $B\in \bk_0$ there are $\delta_n>0$ satisfying $B^{\delta_n}\subseteq U_n$ for infinitely many $n$ (\cite[Proposition 3.1]{dcpdsd}). Let $A_B=\{n:B^{\delta_n}\nsubseteq U_n\}$. Consider the collection $\ac=\{A_B:B\in \bk_0\}$. Let $A_{B_1},\dotsc, A_{B_k}\in \ac$ and on the contrary assume that $\omega\setminus \cup_{i=1}^k A_{B_i}$ is finite. If $B^\prime=\cup_{i=1}^k B_i$ and ${B^\prime}^{\delta_n}\subseteq U_n$, then $B_i^{\delta_n}\subseteq U_n$ and so $n\not\in A_{B_i}$ for all $i=1,\dotsc,k$. Therefore $n\in \omega\setminus \cup_{i=1}^k A_{B_i}$ and hence the set $\{n:{B^\prime}^{\delta_n}\subseteq U_n\}\subseteq \omega\setminus \cup_{i=1}^k A_{B_i}$ is finite which contradicts \cite[Proposition 3.1]{dcpdsd}. Thus our assumption is false and consequently $\ac$ has $\fup$. Clearly $|\ac|<\pf_\square(\jc)$. There exists a $\square$-function $\phi:\omega\rightarrow \omega$ such that $\{n:\phi(n)\in A_B\}\in \jc$ for all $B\in \bk_0$. Proceeding as in $(1)$ then it can be easily shown that $\{U_{\phi(n)}:n\in \omega\}$ is a $\jgbs$-cover of $X$.
\ep

\begin{Th}
\label{TI-2}
Let $\bk$ be a bornology on a metric space $X$ with a closed base $\bk_0$ and let $\ic, \jc$ be ideals on $\omega$. The following statements hold.\\
\noindent$(1)$ If $|\bk_0|<\sle(\ic,\jc)$, then $X$ satisfies $\S1(\IGbs,\JGbs)$.\\
\noindent$(2)$ If $|\bk_0|<\sle(\ast,\jc)$, then $X$ satisfies $\S1(\obs,\JGbs)$, provided $X$ is $\bs$-Lindel\"{o}f.\\
\noindent$(3)$ If $\bk_0|<\slt(\ic,\Fin)$, then $X$ satisfies $\S1(\IGbs,\obs)$.\\
\noindent$(4)$ If $|\bk_0|<\slt(\ast,\Fin)$, then $X$ satisfies $\S1(\obs,\obs)$, provided $X$ is $\bs$-Lindel\"{o}f.
\end{Th}
\bp
$(1)$. Let $\{\uc_n:n\in \omega\}$ be a sequence of $\igbs$-covers of $X$ and let $\uc_n=\{U_{n,m}:m\in \omega\}$. Let $B\in \bk_0$ and $n\in \omega$. Now there exists a sequence $\{\delta_{n,m}:m\in \omega\}$ of positive real numbers such that $\{m:B^{\delta_{n,m}}\nsubseteq U_{n,m}\}\in \ic$. Consider the function $A_B:\omega\rightarrow \pc(\omega)$ defined by $A_B(n)=\{m:B^{\delta_{n,m}}\nsubseteq U_{n,m}\}$. Let $S=\{A_B:B\in \bk_0\}$. Clearly $|S|<\sle(\ic,\jc)$. Choose a $\phi\in {^\omega}\omega$ satisfying $\{n:\phi(n)\in A_B(n)\}\in \jc$ for all $B\in \bk_0$. We claim that $\{U_{n,\phi(n)}:n\in \omega\}$ is a $\jgbs$-cover of $X$. For this we need to show that for $B\in \bk_0$ there is a sequence $\{\sigma_n:n\in \omega\}$ of positive real numbers for which $\{n:B^{\sigma_n}\nsubseteq U_{n,\phi(n)}\}\in \jc$.  Note that $\{n:\phi(n)\in A_B(n)\}\in \jc$. Define $\sigma_n=\delta_{n,\phi(n)}$ for each $n$. Let $T=\{n:B^{\sigma_n}\nsubseteq U_{n,\phi(n)}\}$. For $n\in T$, $B^{\delta_{n,\phi(n)}}\nsubseteq U_{n,\phi(n)}$ and $\phi(n)\in A_B(n)$. Thus $T\subseteq \{n:\phi(n)\in A_B(n)\}$ and hence $T\in \jc$.

$(2)$. Let $\{\uc_n:n\in \omega\}$ be a sequence of open $\bs$-covers of $X$ and let $\uc_n=\{U_{n,m}:m\in \omega\}$. For $B\in \bk_0$ and $n\in \omega$ there are  $\delta_{n,m}>0$ satisfying $B^{\delta_{n,m}}\subseteq U_{n,m}$ for infinitely many $m$. Consider the sequence $A_B:\omega\rightarrow \pc(\omega)$ for which $A_B(n)=\{m:B^{\delta_{n,m}}\nsubseteq U_{n,m}\}$. Let $S=\{A_B:B\in \bk_0\}$. We first show that $S$ has $\fupc$. Take $A_{B_1}, A_{B_2},\dots, A_{B_k}\in S$. Assume that $\omega\setminus \cup_{i=1}^k A_{B_i}(n)=T$ is finite. Let $B^\prime=\cup_{i=1}^k B_i$. Consider the set $\{m:{B^\prime}^{\delta_{n,m}}\subseteq U_{n,m}\}$. Now ${B^\prime}^{\delta_{n,m}}\subseteq U_{n,m}$ implies that $B_i^{\delta_{n,m}}\subseteq U_{n,m}$ for all $i=1,\dotsc, k$. Therefore $m\not\in A_{B_i}(n)$ for all $i=1,\dotsc, k$ and so $m\in T$. Thus $\{m:{B^\prime}^{\delta_{n,m}}\subseteq U_{n,m}\}\subseteq T$ and consequently must be finite which contradicts that $\uc_n$ is an open $\bs$-cover of $X$. Hence $S$ has $\fupc$. Again observe that $|S|<\sle(\ast,\jc)$ and so there exists a $\phi\in {^\omega}\omega$ such that $\{n:\phi(n)\in A_B(n)\}\in \jc$. We now show that $\{U_{n,\phi(n)}:n\in \omega\}$ is a $\jgbs$-cover of $X$. For this let $B\in \bk_0$. Define $\sigma_n=\delta_{n,\phi(n)}$. Clearly $\{n:B^{\sigma_n}\nsubseteq U_{n,\phi(n)}\}\subseteq \{n:\phi(n)\in A_B(n)\}$. Therefore $\{n:B^{\sigma_n}\nsubseteq U_{n,\phi(n)}\}\in \jc$. Hence $X$ satisfies $\S1(\obs,\JGbs)$.

$(3)$. Let $\{\uc_n:n\in \omega\}$ be a sequence of $\igbs$-covers of $X$ and let $\uc_n=\{U_{n,m}:m\in \omega\}$. Proceeding as in $(1)$ we obtain $S=\{A_B:B\in \bk\}$, where $A_B(n)=\{m:B^{\delta_{n,m}}\nsubseteq U_{n,m}\}$. Clearly $|S|<\slt(\ic,\Fin)$. Subsequently there is a $\phi\in {^\omega}\omega$ such that $\{n:\phi(n)\in \cup_{B\in \bk_0}A_B(n)\}\neq \omega$. We show that $\{U_{n.\phi(n)}:n\in \omega\}$ is an open $\bs$-cover of $X$. Let $B\in \bk_0$. There exists a $n$ such that $\phi(n)\not\in \cup_{B\in \bk_0}A_B(n)$. Consequently $B^{\delta_{n,\phi(n)}}\subseteq U_{n,\phi(n)}$. Hence $\{U_{n.\phi(n)}:n\in \omega\}$ is an open $\bs$-cover of $X$.

$(4)$. Let $\{\uc_n:n\in \omega\}$ be a sequence of open $\bs$-covers of $X$ and let $\uc_n=\{U_{n,m}:m\in \omega\}$. Proceeding as in $(2)$ we obtain $S=\{A_B:B\in \bk\}$, where $A_B(n)=\{m:B^{\delta_{n,m}}\nsubseteq U_{n,m}\}$. Clearly $S$ has $\fupc$. As $|S|<\slt(\ast,\Fin)$, choose a $\phi\in {^\omega}\omega$ such that $\{n:\phi(n)\in \cup_{B\in \bk_0}A_B(n)\}\neq \omega$. Then $\{U_{n.\phi(n)}:n\in \omega\}$ is an open $\bs$-cover of $X$.
\ep

\begin{Th}
\label{PI-1}
Let $\bk$ be a bornology on a metric space $X$ with a closed base $\bk_0$ and let $\ic, \jc$ be ideals on $\omega$. The following statements hold.\\
\noindent$(1)$ If $|\bk_0|<\lambda(\ic,\jc)$, then $X$ satisfies $\S1(\IGbs,\JGbs)$.\\
\noindent$(2)$ If $|\bk_0|<\bb_\ic$, then $X$ satisfies $\S1(\Gbs,\IGbs)$.
\end{Th}
\bp
$(1)$. Let $\{\uc_n:n\in \omega\}$ be a sequence of $\igbs$-covers of $X$ and let $\uc_n=\{U_{n,m}:m\in \omega\}$. For $B\in \bk_0$ and $n\in \omega$ there exists a sequence $\{\delta_{n,m}:m\in \omega\}$ of positive real numbers such that $\{m:B^{\delta_{n,m}}\nsubseteq U_{n,m}\}\in \ic$. Let $A_{B,n}=\{m:B^{\delta_{n,m}}\nsubseteq U_{n,m}\}$ and consider the collection $\{A_{B,n}:B\in \bk_0,n\in \omega\}$. Since $|\bk_0|<\lambda(\ic,\jc)$, there exists a $\phi\in {^\omega}\omega$ such that for any $B\in \bk$ $\{n:\phi(n)\in A_{B,n}\}\in \jc$. Clearly $\{U_{n,\phi(n)}:n\in \omega\}$ is a $\jgbs$-cover of $X$ witnessing $\S1(\IGbs,\JGbs)$.

$(2)$. Let $\{\uc_n:n\in \omega\}$ be a sequence of $\gbs$-covers of $X$ and let $\uc_n=\{U_{n,m}:m\in \omega\}$. For $B\in \bk_0$ and $n\in \omega$ there are a $m_0\in \omega$ and a sequence $\{\delta_{n,m}:m\geq m_0\}$ of positive real numbers such that $B^{\delta_{n,m}}\subseteq U_{n,m}$ for all $m\geq m_0$. Choose a function $f_B:\omega\rightarrow \omega$ such that $f_B(n)=\min\{k:\text{ for all } m\geq k \text{ } B^{\delta_{n,m}}\subseteq U_{n,m}\}$. Observe that $|\{f_B:B\in \bk_0\}|<\bb_\ic$. Consequently one can choose a function $g:\omega\rightarrow \omega$ such that $f_B\leq^\ast_\ic g$ for all $B\in \bk_0$. We show that $\{U_{n,g(n)}:n\in \omega\}$ is an $\igbs$-cover of $X$. Let $B\in \bk$. Clearly $f_B\leq^\ast_\ic g$ and so the set $\{n:g(n)\leq f_B(n)\}(=S)\in \ic$. If $n\in S$, then $g(n)\leq f_B(n)$ which implies that $B^{\delta_{n,g(n)}}\nsubseteq U_{n,g(n)}$. If $n\not\in S$, then clearly $B^{\delta_{n,g(n)}}\subseteq U_{n,g(n)}$. Define $\sigma_n=\delta_{n,g(n)}$ for each $n$. It is easy to see that $\{n:B^{\sigma_n}\nsubseteq U_{n,g(n)}\}=S\in \ic$. Therefore $\{U_{n,g(n)}:n\in \omega\}$ is an $\igbs$-cover of $X$. Hence $X$ satisfies $\S1(\Gbs,\IGbs)$.
\ep

\subsection{Results involving ordering on ideals}
In this section we carry out investigation under the orderings $\{1\text{-}1, \KB, \K \}$ on ideals.
\begin{Lemma}
\label{LI}
Let $\bk$ be a bornology on a metric space $X$ with a closed base. Let $\ic_1, \ic_2$ be ideals on $\omega$ with $\ic_1\leq_\psi \ic_2$ for some $\psi\in {^\omega}\omega$. If $\{U_n:n\in \omega\}$ is an $\iogbs$-cover of $X$, then $\{U_{\psi(n)}:n\in \omega\}$ is an $\itgbs$-cover of $X$.
\end{Lemma}
\bp
Let $B\in \bk$. Since $\{U_n:n\in \omega\}$ is an $\iogbs$-cover, there is a sequence $\{\delta_n:n\in \omega\}$ of positive real numbers for which $\{n:B^{\delta_n}\nsubseteq U_n\}\in \ic_1$. Clearly $\psi^{-1}(\{n:B^{\delta_n}\nsubseteq U_n\})(=S)\in \ic_2$ as $\ic_2\leq_\psi \ic_1$. To show that $\{U_{\psi(n)}:n\in \omega\}$ is an $\itgbs$-cover, define $\sigma_n=\delta_{\psi(n)}$ for each $n$. Let $T=\{n:B^{\sigma_n}\nsubseteq U_{\psi(n)}\}$. Note that for $n\in T$, $B^{\delta_{\psi(n)}}\nsubseteq U_{\psi(n)}$ which implies $n\in S$ and so $T\subseteq S$. Therefore $T\in \ic_2$ and hence $\{U_{\psi(n)}:n\in \omega\}$ is an $\itgbs$-cover of $X$.
\ep

The bornological modification of \cite[Proposition 4.4]{js(QN)} is the following.
\begin{Th}
\label{TI}
Let $\bk$ be a bornology on a metric space $X$ with a closed base. Let $\ic_1, \ic_2,\jc_1,\jc_2$ be ideals on $\omega$ with $\ic_1\leq_\K \ic_2$ and $\jc_1\leq_{\KB} \jc_2$. If $X$ satisfies $\S1(\ITGbs,\JOGbs)$, then $X$ satisfies $\S1(\IOGbs,\JTGbs)$.
\end{Th}
\bp
Let $\{\uc_n:n\in \omega\}$ be a sequence of $\iogbs$-covers of $X$ where $\uc_n=\{U_{n,m}:m\in \omega\}$ for $n\in \omega$. Let $\psi:\omega\rightarrow \omega$ be a function satisfying $\ic_1\leq_\psi \ic_2$ and let $\phi:\omega\rightarrow \omega$ be a finite-to-one function satisfying $\jc_1\leq_\phi \jc_2$. We define $V_{n,m}=\cap_{j\in \phi^{-1}(n)}U_{j,\psi(m)}$ if $n\in \phi(\omega)$, and otherwise $V_{n,m}=U_{0,\psi(m)}$. We will show that $\{V_{n,m}:m\in \omega\}$, for each $n$, is an $\itgbs$-cover of $X$. Let $B\in \bk$. Choose a sequence $\{\delta_{n,m}:m\in \omega\}$ of positive real numbers such that $\{m:B^{\delta_{n,m}}\nsubseteq U_{n,m}\}\in \ic_1$ and $\psi^{-1}(\{m:B^{\delta_{n,m}}\nsubseteq U_{n,m}\})\in \ic_2$ for each $n\in \omega$. Now define $\sigma_{n,m}=\min\{\delta_{j,\psi(m)}:j\in \phi^{-1}(n)\}$ for $n,m\in \omega$. It easy to see observe that $\{m:B^{\delta_{n,m}}\nsubseteq V_{n,m}\}\subseteq \cup_{j\in \phi^{-1}(n)}\{m:B^{\delta_{j,\psi(m)}}\nsubseteq U_{j,\psi(m)}\}$ and consequently $\{m:B^{\sigma_{n,m}}\nsubseteq V_{n,m}\}\subseteq \cup_{j\in \phi^{-1}(n)}\psi^{-1}(\{m:B^{\delta_{j,m}}\nsubseteq U_{j,m}\}$. Therefore $\{m:B^{\sigma_{n,m}}\nsubseteq V_{n,m}\}\in \ic_2$ and so  $\{V_{n,m}:m\in \omega\}$ is an $\itgbs$-cover of $X$ for each $n\in \omega$. Applying $\S1(\ITGbs,\JOGbs)$ to $\{V_{n,m}:m,n\in \omega\}$, we can obtain a sequence $\{k_n:n\in \omega\}$ such that $\{V_{n,k_n}:n\in \omega\}$ would be a $\iogbs$-cover of $X$. Define $m_n=\psi(k_{\phi(n)})$ for each $n$. We claim that $\{U_{n,{m_n}}:n\in \omega\}$ is a $\jtgbs$-cover of $X$. Note that for $B\in \bk$ there is a sequence $\{\varepsilon_n:n\in \omega\}$ of positive real numbers such that $\{n:B^{\varepsilon_n}\nsubseteq V_{n,{k_n}}\}\in \jc_1$. Clearly $\phi^{-1}(\{n:B^{\varepsilon_n}\nsubseteq V_{n,{k_n}}\})\in \jc_2$. Define $\eta_n=\varepsilon_{\phi(n)}$ for each $n$. Now it is easy to verify that $\{n:B^{\eta_n}\nsubseteq U_{n,{m_n}}\}\subseteq \{n:B^{\varepsilon_{\phi(n)}}\nsubseteq V_{\phi(n),k_{\phi(n)}}\}$. Consequently $\{n:B^{\eta_n}\nsubseteq U_{n,{m_n}}\}\subseteq \phi^{-1}(\{n:B^{\varepsilon_n}\nsubseteq V_{n,{k_n}}\})$ and so $\{n:B^{\eta_n}\nsubseteq U_{n,{m_n}}\}\in \jc_2$. This shows that $\{U_{n,{m_n}}:n\in \omega\}$ is a $\jtgbs$-cover of $X$ witnessing $\S1(\IOGbs,\JTGbs)$.
\ep

\begin{Th}
\label{TI-}
Let $\bk$ be a bornology on a metric space $X$ with a closed base. Let $\ic_1, \ic_2,\jc_1,\jc_2$ be ideals on $\omega$ with $\ic_1\leq_\K \ic_2$ and $\jc_1\leq_{\KB} \jc_2$. The following statements hold.\\
\noindent$(1)$ If $X$ satisfies $\S1(\ITGbs,\obs)$, then $X$ satisfies $\S1(\IOGbs,\obs)$.\\
\noindent$(2)$ If $X$ satisfies $\S1(\obs,\JOGbs)$, then $X$ satisfies $\S1(\obs,\JTGbs)$.
\end{Th}
\bp
$(1)$. Let $\psi:\omega\rightarrow \omega$ be a function satisfying $\ic_1\leq_\psi \ic_2$. Start with $\{\uc_n:n\in \omega\}$, a sequence of $\iogbs$-covers of $X$ where $\uc_n=\{U_{n,m}:m\in \omega\}$. By Lemma \ref{LI}, $\vc_n=\{U_{n,\psi(m)}:m\in \omega\}$ is an $\itgbs$-cover of $X$. Applying $\S1(\ITGbs,\obs)$ to $\{\vc_n:n\in \omega\}$, we obtain a sequence $\{m_n:n\in \omega\}$ such that $\{U_{n,\phi(m_n)}:n\in \omega\}$ is an open $\bs$-cover of $X$. Hence $X$ satisfies $\S1(\IOGbs,\obs)$.

$(2)$. Let $\phi$ be a finite-to-one function witnessing $\jc_1\leq_\phi \jc_2$. Let $\{\uc_n:n\in \omega\}$ be a sequence of open $\bs$-covers of $X$ and let $\uc_n=\{U_{n,m}:m\in \omega\}$. If $i\in \phi(\omega)$ and $i=\phi(n)$, take $\vc_i=\{U_{1,m_1}\cap \dotsc \cap U_{n,m_n}:U_{1,m_1}\in \uc_1,\dotsc, U_{n,m_n}\in \uc_n\}$ and take $\vc_i=\{U_{0,m}:m\in \omega\}$ otherwise. Applying $\S1(\obs,\JOGbs)$ to $\{\vc_i:i\in \omega\}$, one can find a sequence $\{k_i:i\in \omega\}$ such that $\{V_{i,k_i}:i\in \omega\}$ is a $\jogbs$-cover. Now observe that there is a sequence $\{m_n:n\in \omega\}$ such that $V_{\phi(n),k_{\phi(n)}}\subseteq U_{n,m_n}$. We now show that $\{U_{n,m_n}:n\in \omega\}$ is a $\jtgbs$-cover. Let $B\in \bk$. Choose a sequence $\{\delta_i:i\in \omega\}$ of positive real numbers such that $\{i:B^{\delta_i}\nsubseteq V_{i,k_i}\}\in \jc_1$. Clearly $\phi^{-1}(\{i:B^{\delta_i}\nsubseteq V_{i,k_i}\})\in \jc_2$. Define $\sigma_n=\delta_{\phi(n)}$ for each $n\in \omega$. It is now easy to check that $\{n:B^{\sigma_n}\nsubseteq U_{n,m_n}\}\subseteq \{n:B^{\delta_{\phi(n)}}\nsubseteq V_{\phi(n),k_{\phi(n)}}\}$. Consequently $\{n:B^{\sigma_n}\nsubseteq U_{n,m_n}\}\subseteq \phi^{-1}(\{i:B^{\delta_i}\nsubseteq V_{i,k_i}\})\in \jc_2$. Therefore $\{U_{n,m_n}:n\in \omega\}$ is a $\jtgbs$-cover. Hence $X$ satisfies $\S1(\obs,\JTGbs)$.
\ep

When $\ic_1\leq_\K \ic_2$ and $\jc_1\leq_{\KB} \jc_2$, the relations among the selection principles $\S1(\ac,\bc)$ obtained above are represented in Figure \ref{diag3}.
\begin{figure}[h]
 \begin{adjustbox}{keepaspectratio,center}
\begin{tikzcd}[column sep=2ex,row sep=5ex,arrows={crossing over}]
\S1(\obs,\JTGbs)\arrow[r]&\S1(\IOGbs,\JTGbs)\arrow[r]&\S1(\IOGbs,\obs)&\\
\S1(\obs,\JOGbs)\arrow[r]\arrow[u]&\S1(\ITGbs,\JOGbs)\arrow[r]\arrow[u]&\S1(\IOGbs,\obs)\arrow[u]
\end{tikzcd}
\end{adjustbox}
 \caption{Diagram of the selection principles when $\ic_1\leq_\K \ic_2$ and $\jc_1\leq_{\KB} \jc_2$}
 \label{diag3}
\end{figure}

\begin{Th}
\label{TI-}
Let $\bk$ be a bornology on a metric space $X$ with a closed base. Let $\ic_1, \ic_2,\jc_1,\jc_2$ be ideals on $\omega$ with $\ic_1\leq_\square \ic_2$ and $\jc_1\leq_\square \jc_2$, where $\square\in \{1\text{-}1, \KB, \K \}$. The following statements hold.\\
\noindent$(1)$ If $X$ satisfies $[\ITGbs,\JOGbs]_\square$, then $X$ satisfies $[\IOGbs,\JTGbs]_\square$.\\
\noindent$(2)$ If $X$ satisfies $[\obs,\JOGbs]_\square$, then $X$ satisfies $[\obs,\JTGbs]_\square$, provided $X$ is $\bs$-Lindel\"{o}f.
\end{Th}
\bp
We prove only $(1)$. Let $\ic_1\leq_\psi \ic_2$ and $\jc_1\leq_\phi \jc_2$, where $\psi, \phi$ are $\square$-functions. Let $\{U_n:n\in \omega\}$ be an $\iogbs$-cover of $X$. By Lemma \ref{LI}, $\{U_{\psi(n)}:n\in \omega\}$ is an $\itgbs$-cover of $X$. By the given condition, there exists a function $\alpha:\omega\rightarrow \omega$ such that $\{U_{\psi(\alpha(n))}:n\in \omega\}$ is a $\jogbs$-cover. Again by Lemma \ref{LI}, $\{U_{\psi(\alpha(\phi(n))))}:n\in \omega\}$ is a $\jtgbs$-cover of $X$. Hence $X$ satisfies $[\IOGbs,\JTGbs]_\square$.
\ep

Theorem \ref{TI-} can be represented by Figure \ref{diag2} where the ideals $\ic_1,\ic_2,\jc_1,\jc_2$ with $\ic_1\leq_\square \ic_2$ and $\jc_1\leq_\square \jc_2$, where $\square\in \{1-1,\KB,\K\}$.
\begin{figure}[h]
 \begin{adjustbox}{keepaspectratio,center}
\begin{tikzcd}
\left[\obs,\JTGbs\right]_\square\arrow[r] & \left[\ITGbs,\JTGbs\right]_\square\\
\left[\obs,\JOGbs\right]_\square\arrow[r]\arrow[u] & \left[\ITGbs,\JOGbs\right]_\square\arrow[u]
\end{tikzcd}
\end{adjustbox}
 \caption{Diagram of the selection principles $[\ac,\bc]_\square$}
 \label{diag2}
\end{figure}

\begin{Th}
\label{TI-}
Let $\bk$ be a bornology on a metric space $X$ with a closed base. Let $\ic_1, \ic_2$ be ideals on $\omega$ with $\ic_1\leq_\K \ic_2$. The following statements hold.\\
\noindent$(1)$ $X$ satisfies $[\obs,\ITGbs]_\K$ and $\S1(\ITGbs,\IOGbs)$ if and only if $X$ satisfies $\S1(\obs, \IOGbs)$.\\
\noindent$(2)$ $X$ satisfies $[\ITGbs,\IOGbs]_\K$ and $\S1(\IOGbs,\IOGbs)$ if and only if $X$ satisfies $\S1(\ITGbs, \IOGbs)$.
\end{Th}
\bp
We prove only $(2)$. Let $\{\uc_n:n\in \omega\}$ be a sequence of $\itgbs$-covers of $X$ and write $\uc_n=\{U_{n,m}:m\in \omega\}$ for $n\in \omega$. Clearly for each $n$ there exists a function $\alpha_n:\omega\rightarrow \omega$ such that $\{U_{n,\alpha_n(m)}:m\in \omega\}(=\vc_n)$ is an $\iogbs$-cover of $X$. Now applying $\S1(\IOGbs,\IOGbs)$ to $\{\vc_n:n\in \omega\}$, we can find a sequence $\{m_n:n\in \omega\}$ such that $\{U_{n,\alpha_n(m_n)}:n\in \omega\}$ is an $\iogbs$-cover of $X$. Hence $X$ satisfies $\S1(\ITGbs, \IOGbs)$.

Conversely, let $\ic_1\leq_\psi \ic_2$ for some $\psi\in {^\omega}\omega$ and let $\{\uc_n:n\in \omega\}$ be a sequence of $\iogbs$-covers of $X$ where $\uc_n=\{U_{n,m}:m\in \omega\}$. As $\ic_1\leq_\psi \ic_2$, $\vc_n=\{U_{n,\psi(m)}:m\in \omega\}$ would be an $\itgbs$-cover. Applying $\S1(\ITGbs, \IOGbs)$ to $\{\vc_n:n\in \omega\}$, we can choose a sequence $\{m_n:n\in \omega\}$ such that $\{U_{n,\psi(m_n)}:n\in \omega\}$ is an $\iogbs$-cover of $X$. Hence $X$ satisfies $\S1(\IOGbs,\IOGbs)$. Also it is easy to show that $X$ satisfies $[\ITGbs,\IOGbs]_\K$.
\ep

Now consider an ideal $\ic$ which is not tall. Let $\square\in \{1\text{-}1, \KB, \K \}$. In \cite[Remark 4.2]{js(QN)}, \v{S}upina observed that $\ic$ is not tall if and only if $\ic\leq_\square \Fin$ and we use this fact in the following lemma.
\begin{Lemma}
\label{LI-2}
Let $\bk$ be a bornology on a metric space $X$ with a closed base and let $\ic$ be not a tall ideal on $\omega$. If $\{U_n:n\in \omega\}$ is an $\igbs$-cover of $X$, then $\{U_{\phi(n)}:n\in \omega\}$ is a $\gbs$-cover of $X$ for some $\phi\in {^\omega}\omega$.
\end{Lemma}
\bp
As $\ic$ is not a tall ideal, $\ic\leq_\square \Fin$. There is a $\square$-function $\phi$ witnessing $\ic\leq_\phi \Fin$. Since $\{U_n:n\in \omega\}$ is an $\igbs$-cover, $\{U_{\phi(n)}:n\in \omega\}$ is a $\Fin$-$\gbs$-cover by Lemma \ref{LI}. Hence $\{U_{\phi(n)}:n\in \omega\}$ is a $\gbs$-cover of $X$.
\ep

The next result can be easily verified.
\begin{Th}
\label{TI-}
Let $\bk$ be a bornology on a metric space $X$ with a closed base. Let $\ic, \jc$ be tall ideals on $\omega$ and $\square\in \{1\text{-}1, \KB, \K \}$. The following statements hold.\\
\noindent$(1)$ $X$ satisfies $\S1(\obs, \JGbs)$ if and only if $X$ satisfies $\S1(\obs,\Gbs)$, provided $X$ is $\bs$-Lindel\"{o}f.\\
\noindent$(2)$ $X$ satisfies $\S1(\IGbs,\obs)$ if and only if $X$ satisfies $\S1(\Gbs, \obs)$.\\
\noindent$(3)$ $X$ satisfies $[\IGbs, \JGbs]_\square$ if and only if $X$ satisfies $[\IGbs,\Gbs]_\square$.\\
\noindent$(4)$ $X$ satisfies $[\obs, \JGbs]_\square$ if and only if $X$ satisfies $[\obs,\Gbs]_\square$, provided $X$ is $\bs$-Lindel\"{o}f.
\end{Th}

Also by Lemma \ref{LI-2}, we have the following.
\begin{Prop}
\label{PI-1}
Let $\bk$ be a bornology on a metric space $X$ with a closed base.\\
\noindent$(1)$ Let $\ic_1,\ic_2$ be ideals on $\omega$ with $\ic_2\subseteq \ic_1$ and $\jc_2\subseteq \jc_1$. If $X$ satisfies $\S1(\IOGbs,\JTGbs)$, then $X$ satisfies $\S1(\ITGbs,\JOGbs)$.\\
\noindent$(2)$ If $X$ satisfies $\S1(\obs,\Gbs)$, then $X$ satisfies $\S1(\IGbs,\JGbs)$ for any ideals $\ic,\jc$.
\end{Prop}

\subsection{Results on selection principles involving ideals having a pseudounion}
The notion of pseudounion of a family $\ac$ has been introduced and studied in \cite{lbpdjs}. Suppose that $\ic,\jc$ are ideals and $\jc$ has a pseudounion, then $X$ satisfies $\S1(\ic\text{-}\Gamma,\jc\text{-}\Gamma)$ if and only if $\S1(\ic\text{-}\Gamma,\Gamma)$ \cite[Theorem 5.3]{lbpdjs}. For our bornological investigation, in this section we consider ideals having a pseudounion. The following lemma will be used subsequently.
\begin{Lemma}
\label{LIP-1}
Let $\bk$ be a bornology on a metric space $X$ with a closed base and let $\ic$ be an ideal on $\omega$. If $\ic$ has a pseudounion, then every $\igbs$-cover of $X$ has a subset which is a $\gbs$-cover of $X$.
\end{Lemma}
\bp
Let $A$ be a pseudounion of $\ic$ and let $C=\omega\setminus A$. Since both $A$ and $C$ are infinite, we can find an injection $\sigma:A\rightarrow C$.  Let $\uc=\{U_n:n\in \omega\}$ be an $\igbs$-cover of $X$. If $n\in C$, choose $m_n=n$, and if $n\in A$, choose $m_n=\sigma(n)$. We claim that $\{U_{m_n}:n\in \omega\}$ is a $\gbs$-cover of $X$. Let $B\in \bk$. Then there is a sequence $\{\delta_n:n\in \omega\}$ of positive real numbers such that $\{n:B^{\delta_n}\nsubseteq U_n\}(=K)\in \ic$. Define $\varepsilon_n=\delta_n$ if $n\in C$ and $\varepsilon_n=\delta_{\sigma(n)}$ if $n\in A$. Now consider the set $L=\{n:B^{\varepsilon_n}\nsubseteq U_{m_n}\}$. We have $K\subseteq^* A$ which implies $K\setminus A$ is finite and consequently, $C\cap K$ is finite and so $C\cap L=C\cap K$. If $n\in A$, then $B^{\varepsilon_n}\nsubseteq U_{m_n}$ implies that $B^{\delta_{\sigma(n)}}\nsubseteq U_{\sigma(n)}$. Consequently, $\sigma(n)\in K$ and so $\sigma(n)\in C\cap K$. Therefore for any $n\in A\cap L$ $\sigma(n)\in C\cap K$. Since $\sigma$ is injection and $C\cap K$ is finite, $A\cap L$ must also be finite. As we already have $C\cap L$ is finite, hence $L$ must be finite which shows that $\{U_{m_n}:n\in \omega\}$ is a $\gbs$-cover of $X$.
\ep

\begin{Th}
\label{TIP-1}
Let $\bk$ be a bornology on a metric space $X$ with a closed base. Let $\ic, \jc$ be ideals on $\omega$ and $\jc$ have a pseudounion. The following statements are equivalent.\\
\noindent$(1)$ $X$ satisfies $\S1(\IGbs,\JGbs)$.\\
\noindent$(2)$ $X$ satisfies $\S1(\IGbs,\Gbs)$.
\end{Th}
\bp
We prove only $(1)\Rightarrow (2)$. Let $A$ be a pseudounion of $\jc$ and let $C=\omega\setminus A$. Since $A$ and $C$ are infinite, we can choose an injection $\sigma:A\rightarrow C$ such that $\sigma(n)>n$ for each $n\in A$. Let $\{\uc_n:n\in \omega\}$ be a sequence of $\igbs$-covers of $X$ and let $\uc_n=\{U_{n,m}:m\in \omega\}$. Construct a sequence $\{\vc_n:n\in \omega\}$ of $\igbs$-cover of $X$ as follows. Define $\vc_1=\uc_1$ and $\vc_{n+1}=\{V_{n+1,m}:m\in \omega\}$, where $V_{n+1,m}=V_{n,m}\cap U_{n+1,m}$. Applying $\S1(\IGbs,\JGbs)$ to $\{\vc_n:n\in \omega\}$, we can obtain a $m_n$ for each $n$ such that $\{V_{n,m_n}:n\in \omega\}$ is a $\jgbs$-cover of $X$. Define $k_n=m_n$ if $n\in C$ and $k_n=m_{\sigma(n)}$ if $n\in A$. We show that $\{V_{n,k_n}:n\in \omega\}$ is a $\gbs$-cover of $X$.Take a $B\in \bk$ and choose a sequence $\{\delta_n:n\in \omega\}$ of positive real numbers such that $\{n:B^{\delta_n}\nsubseteq V_{n,m_n}\}(=K)\in \ic$. Define $\varepsilon_n=\delta_n$ if $n\in C$ and $\varepsilon_n=\delta_{\sigma(n)}$ if $n\in A$. Let $L=\{n:B^{\varepsilon_n}\nsubseteq V_{n,k_n}\}$. Proceeding as in Lemma \ref{LIP-1}, we can show that $C\cap L$ is finite. If $n\in A$, $V_{\sigma(n),m_{\sigma(n)}}\subseteq V_{n,m_{\sigma(n)}}=V_{n,k_n}$ as $\sigma(n)>n$. Consequently $B^{\varepsilon_n}\nsubseteq V_{n,k_n}$ implies that $B^{\varepsilon_n}\nsubseteq V_{\sigma(n),m_{\sigma(n)}}$. Therefore $\sigma(n)\in K$. Clearly if $n\in A\cap L$, then $\sigma(n)\in C\cap K$. Hence $A\cap L$ must also be finite. Since $C\cap L$ is already finite, $L$ is finite. This shows that $\{V_{n,k_n}:n\in \omega\}$ is a $\gbs$-cover of $X$. As $V_{n,k_n}=V_{n-1,k_n}\cap U_{n,k_n}$ for each $n$, $\{U_{n,k_n}:n\in \omega\}$ forms the required $\gbs$-cover of $X$. Hence $X$ satisfies $\S1(\IGbs,\Gbs)$.
\ep

\begin{Th}
\label{TIP-3}
Let $\bk$ be a bornology on a metric space $X$ with a closed base and let $X$ be $\bs$-Lindel\"{o}f. Let the ideal $\jc$ on $\omega$ have a pseudounion. The following statements are equivalent.\\
\noindent$(1)$ $X$ satisfies $\S1(\obs,\JGbs)$.\\
\noindent$(2)$ $X$ satisfies $\S1(\obs,\Gbs)$.
\end{Th}
\bp
We prove $(1)\Rightarrow (2)$. Let $\{\uc_n:n\in \omega\}$ be a sequence of open $\bs$-covers of $X$ and $\uc_n=\{U_{n,m}:m\in \omega\}$. Take $\vc_n=\{U_{1,m_1}\cap \dotsc \cap U_{n,m_n}:U_{1,m_1}\in \uc_1,\dotsc,U_{n,m_n}\in \uc_n\}$ for $n\in \omega$. By \cite[Lemma 3.1]{dcpdsd}, $\vc_n$ is an open $\bs$-cover of $X$ for each $n$. Applying $\S1(\obs,\JGbs)$ to $\{\vc_n:n\in \omega\}$, we can choose a $V_n\in \vc_n$ such that $\{V_n:n\in \omega\}$ forms a $\jgbs$-cover of $X$. Now using Lemma \ref{LIP-1}, we can find a subsequence $\{V_{n_k}:k\in \omega\}$ which is a $\gbs$-cover of $X$. From the construction of $V_{n_k}$ we can then choose $U_{i,m_i}$ where $i\leq n_1$ and $U_{i,m_i}$ where $n_{k-1}<i\leq n_k$, $k>1$. It is easy to verify that $\{U_{i,m_i}:i\in \omega\}$ is a $\gbs$-cover of $X$. Hence $(2)$ holds.
\ep

\begin{Th}
\label{TIP-2}
Let $\bk$ be a bornology on a metric space $X$ with a closed base. Let the ideals $\ic, \jc$ on $\omega$ have a pseudounion. The following statements are equivalent.\\
\noindent$(1)$ $X$ satisfies $\S1(\IGbs,\JGbs)$.\\
\noindent$(2)$ $X$ satisfies $\S1(\Gbs,\Gbs)$.
\end{Th}
\bp
$(1)\Rightarrow (2)$. By Theorem \ref{TIP-1}, $X$ satisfies $\S1(\IGbs,\Gbs)$. Clearly $X$ satisfies $\S1(\Gbs,\Gbs)$.

$(2)\Rightarrow (1)$. In view of Lemma \ref{LIP-1}, it follows that $X$ satisfies $\S1(\IGbs,\Gbs)$. Now by Theorem \ref{TIP-1}, $X$ satisfies $\S1(\IGbs,\JGbs)$.
\ep

Using the same line of argument we can prove the following.

\begin{Th}
\label{TIP-4}
Let $\bk$ be a bornology on a metric space $X$ with a closed base. Let the ideal $\jc$ on $\omega$ have a pseudounion. The following statements are equivalent.\\
\noindent$(1)$ $X$ satisfies $\S1(\Gbs,\JGbs)$.\\
\noindent$(2)$ $X$ satisfies $\S1(\Gbs,\Gbs)$.
\end{Th}

\begin{Th}
\label{TIP-5}
Let $\bk$ be a bornology on a metric space $X$ with a closed base. Let the ideal $\ic$ on $\omega$ have a pseudounion. The following statements are equivalent.\\
\noindent$(1)$ $X$ satisfies $\S1(\IGbs,\Gbs)$.\\
\noindent$(2)$ $X$ satisfies $\S1(\Gbs,\Gbs)$.
\end{Th}

\begin{Th}
\label{TIP-6}
Let $\bk$ be a bornology on a metric space $X$ with a closed base. Let the ideal $\ic$ on $\omega$ have a pseudounion. The following statements are equivalent.\\
\noindent$(1)$ $X$ satisfies $\S1(\IGbs,\obs)$.\\
\noindent$(2)$ $X$ satisfies $\S1(\Gbs,\obs)$.
\end{Th}

The above implications among selection principles are describe in the following diagram (Figure \ref{diag4}).
\begin{figure}[h]
 \begin{adjustbox}{keepaspectratio,center}
\begin{tikzcd}[column sep=2ex,row sep=5ex,arrows={crossing over}]
\S1(\obs,\JGbs)\arrow[r]&\S1(\IGbs,\JGbs)\arrow[r,thick,equal]&\S1(\IGbs,\Gbs)\arrow[r]&\S1(\IGbs,\obs)&\\
\S1(\obs,\Gbs)\arrow[r]\arrow[u,thick,equal]&\S1(\Gbs,\JGbs)\arrow[r,thick,equal]\arrow[u,thick,equal]&\S1(\Gbs,\Gbs)\arrow[r]\arrow[u,thick,equal]&\S1(\Gbs,\obs)\arrow[u,thick,equal]&\\
\end{tikzcd}
\end{adjustbox}
 \caption{Diagram of the selection principles when $\ic,\jc$ have pseudounion}
 \label{diag4}
\end{figure}

We conclude this section with certain observations on splittability. Every spaces satisfy $\Split(\Gbs,\Gbs)$ and $\Split(\Gbs,\JGbs)$ for any ideal $\jc$.
If $\ic$ has a pseudounion, then every $\igbs$-cover of $X$ has a subset which is a $\gbs$-cover (Lemma \ref{LIP-1}). Therefore every spaces satisfy $\Split(\IGbs,\JGbs)$ and $\Split(\IGbs,\Gbs)$.

Combining Theorem \ref{TIP-3} and \cite[Proposition 2.4]{dcpdsd2}, we obtain the following result.
\begin{Prop}
\label{PIS-1}
Let $\bk$ be a bornology on a metric space $X$ with a closed base and let the ideal $\jc$ on $\omega$ have a pseudounion. The following statements are equivalent.\\
\noindent$(1)$ $X$ satisfies $\S1(\obs,\JGbs)$.\\
\noindent$(2)$ $X$ satisfies $\S1(\obs,\Gbs)$.\\
\noindent$(3)$ $X$ satisfies $\Split(\obs,\Gbs)$.\\
\noindent$(4)$ $X$ satisfies $\Split(\obs,\JGbs)$.
\end{Prop}

For any two ideals $\ic$ and $\jc$ the implications among $\Split(\ac,\bc)$, where $\ac,\bc\in \{\obs,\Gbs,\IGbs,\JGbs\}$ are describe in Figure \ref{diag7}.

\begin{figure}[h]
 \begin{adjustbox}{keepaspectratio,center}
\begin{tikzcd}[column sep=2ex,row sep=5ex,arrows={crossing over}]
&\Split(\IGbs,\obs)\\
\Split(\obs,\JGbs)\arrow[r]&\Split(\IGbs,\JGbs)\arrow[r]\arrow[u]&\Split(\Gbs,\JGbs)\\
&\Split(\IGbs,\Gbs)\arrow[r]\arrow[u]&\Split(\Gbs,\Gbs)\arrow[u,thick,equals]
\end{tikzcd}
\end{adjustbox}
 \caption{Diagram of splittability with respect to ideals $\ic,\jc$}
 \label{diag7}
\end{figure}

\section{The $\ic$-$\bs$-Hurewicz Property}
The $\ic$-$\bs$-Hurewicz Property of $X$ introduced in \cite{sd}. As usual we consider ideals which have a pseudounion and orders $\{1\text{-}1,\KB,\K\}$ on ideals. Next two propositions are modifications of \cite[Propositions 5.1, 5.2]{dcpdsd2}.
\begin{Prop}
\label{PIH-1}
Let $\bk$ be a bornology on a metric space $X$ with a closed base and $Y$ be another metric space and let $\ic$ be an ideal on $\omega$. Let $f:X\rightarrow Y$ be a continuous function. If $X$ has the $\bs$-Hurewicz property, then $f(X)$ has the $\ibs$-Hurewicz property.
\end{Prop}

\begin{Prop}
\label{PIH-2}
Let $\bk$ be a bornology on a metric space $X$ with a closed base and let $\ic$ be an ideal on $\omega$. If $X$ has the $\ibs$-Hurewicz property, then every continuous image of $X$ into $\omega^\omega$ is $\ic$-bounded.
\end{Prop}

\begin{Prop}
\label{PIH-3}
Let $\bk$ be a bornology on a metric space $X$ with a closed base $\bk_0$ and let $X$ be $\bs$-Lindel\"{o}f. Let $\ic$ be an ideal on $\omega$. If $|\bk_0|<\bb_\ic$, then $X$ has the $\ibs$-Hurewicz property.
\end{Prop}
\bp
Let $\{\uc_n:n\in \omega\}$ be a sequence of open $\bs$-covers of $X$, where $\uc_n=\{U_{n,m}:m\in \omega\}$ for each $n$.	 For $B\in \bk_0$ let $f_B(n)=\min\{m:B^\delta\subseteq U_{n,m} \text{ for some } \delta>0\}$. Consider the set $\{f_B:B\in \bk_0\}$. As $|\bk_0|<\bb_\ic$ choose a function $g:\omega\rightarrow \omega$ satisfying $f_B\leq^\ast_\ic g$. Let $\vc_n=\{U_{n,m}:m\leq g(n)\}$. Each $\vc_n$ is a finite subset of $\uc_n$. Let $B\in \bk_0$. We have $f_B\leq^\ast_\ic g$ and $S=\{n:g(n)<f_B(n)\}\in \ic$. If $n\in S$, then $g(n)<f_B(n)$ and so $B^\delta\nsubseteq U_{n,m}$ for any $\delta>0$. Choose a $\delta_n>0$ with $B^{\delta_n}\nsubseteq U_{n,m}$. When $n\not\in S$, $f_B(n)\leq g(n)$ and consequently $B^{\delta_{f_B(n)}}\subseteq U_{n,f_B(n)}$. Define $\sigma_n=\delta_n$ if $n\in S$ and $\sigma_n=\delta_{f_{B_0}(n)}$ if $n\not\in S$. Now it is easy to verify that $\{n:B^{\delta_n}\nsubseteq U_{n,m} \text{ for any } U_{n,m}\in \vc_n\}=S\in \ic$. Hence $\{\vc_n:n\in \omega\}$ witnesses the $\ibs$-Hurewicz property of $X$.
\ep

\begin{Th}
\label{TIH-1}
Let $\bk$ be a bornology on a metric space $X$ with a closed base. Let $\ic_1, \ic_2$ be ideals on $\omega$ with $\ic_1\leq_\square \ic_2$, where $\square\in \{1\text{-}1,\KB,\K\}$ .
If $X$ has the $\iobs$-Hurewicz property, then $X$ has the $\itbs$-Hurewicz property.
\end{Th}
\bp
Let $\ic_1\leq_\phi \ic_2$ for some $\phi\in {^\omega}\omega$. Let $\{\uc_n:n\in \omega\}$ be a sequence of open $\bs$-covers of $X$, where $\uc_n=\{U_{n,m}:m\in \omega\}$ for each $n$. If $i\in \phi(\omega)$ and $i=\phi(n)$, let $\vc_i=\{U_{1,m_1}\cap\dotsc \cap U_{n,m_n}:U_{1,m_1}\in \uc_1,\dotsc, U_{n,m_n}\in \uc_n\}$ otherwise let $\vc_i=\{U_{0,m}:m\in \omega\}$. Since $X$ has the $\iobs$-Hurewicz property, choose a sequence $\{\wc_n:n\in \omega\}$ of finite sets with $\wc_n\subseteq \vc_n$ for each $n$ such that for $B\in \bk$ there sequence $\{\delta_n:n\in \omega\}$ of positive real numbers satisfying $\{n:B^{\delta_n}\nsubseteq U \text{ for any } U\in \vc_n\}\in \ic_1$. Let $\zc_n$ be the collection of all the $n$th components of members of $\wc_{\phi(n)}$. Clearly $\zc_n$ is a finite subset of $\uc_n$. We now show that $\{\zc_n:n\in \omega\}$ witnesses the $\itbs$-Hurewicz property. Let $B\in \bk$. As $\ic_1\leq_\phi \ic_2$, we have $\phi^{-1}(\{n:B^{\delta_n}\nsubseteq U \text{ for any } U\in \vc_n\})\in \ic_2$. Define $\sigma_n=\delta_{\phi(n)}$ for each $n$. The proof will be complete if we show that $\{n:B^{\sigma_n}\nsubseteq U \text{ for any } U\in \zc_n\}\in \ic_2$. It is easy to see that for each $U\in \zc_n$ there is a $V\in \wc_{\phi(n)}$ satisfying $V\subseteq U$. Therefore $\{n:B^{\delta_n}\nsubseteq U \text{ for any } U\in \zc_n\}\subseteq\{n:B^{\delta_{\phi(n)}}\nsubseteq V \text{ for any } V\in \wc_{\phi(n)}\}$. Consequently $\{n:B^{\delta_n}\nsubseteq U \text{ for any } U\in \zc_n\}\subseteq \phi^{-1}(\{n:B^{\delta_n}\nsubseteq U \text{ for any } U\in \wc_n\}\in \ic_2$.
\ep

\begin{Th}
\label{TIH-2}
Let $\bk$ be a bornology on a metric space $X$ with a closed base. Let $\ic_1, \ic_2$ be ideals on $\omega$ with $\ic_1\leq_\square \ic_2$, where $\square\in \{1\text{-}1,\KB,\K\}$ .
If ONE has no winning strategy in the $\iobs$-Hurewicz game, then ONE has no winning strategy in the $\itbs$-Hurewicz game.
\end{Th}
\bp
Let $\ic_1\leq_\phi \ic_2$ for some $\phi\in {^\omega}\omega$. Let $F$ be a strategy for ONE in the $\itbs$-Hurewicz game. We define a strategy $G$ for ONE in the
$\iobs$-Hurewicz game as follows. Let $F(\emptyset)=\{U_{(n)}:n\in \omega\}$, an open $\bs$-cover of $X$. Define $G(\emptyset)=\{U_{(n)}:n\in \omega\}$. Let TWO choose a finite subset $\vc_{(n_1)}=\{U_{(n)}:n\leq n_1\}$ of $G(\emptyset)$. In the $\itbs$-Hurewicz game, let TWO choose $\vc_{n_{\phi(1)}}=\{U_{(n)}:n\leq n_{\phi(1)}\}$, a finite subset of $F(\emptyset)$. Suppose that $U_{(n_{\phi(1)},\dotsc, n_{\phi(m-1)})}$ is chosen. Let $F(U_{(n_{\phi(1)},\dotsc, n_{\phi(m-1)})})=\{U_{(n_{\phi(1)},\dotsc, n_{\phi(m-1)},n)}:n\in \omega\}$. Define $G(U_{(n_{\phi(1)},\dotsc, n_{\phi(m-2)},n_{m-1})})=\{U_{(n_{\phi(1)},\dotsc, n_{\phi(m-1)},n)}:n\in \omega\}$. Let TWO choose a finite set $\vc_{(n_{\phi(1)},\dotsc, n_{\phi(m-1)},n_m)}=\{U_{(n_{\phi(1)},\dotsc, n_{\phi(m-1)},n)}:n\leq n_m\}$. In the $\itbs$-Hurewicz game, let TWO choose $\vc_{(n_{\phi(1)},\dotsc ,n_{\phi(m)})}=\{U_{(n_{\phi(1)},\dotsc, n_{\phi(m-1)},n)}:n\leq n_{\phi(m)}\}$ and so on. This defines the strategy $G$. A play in the $\itbs$-Hurewicz game is $$F(\emptyset),\vc_{n_{\phi(1)}},\dotsc, F(U_{(n_{\phi(1)},\dotsc, n_{\phi(m-1)})}), \vc_{(n_{\phi(1)},\dotsc ,n_{\phi(m)})},\dotsc$$ Correspondingly a play in the $\iobs$-Hurewicz game is $$G(\emptyset),\vc_{(n_1)},\dotsc, G(U_{(n_{\phi(1)},\dotsc, n_{\phi(m-2)},n_{m-1})}), \vc_{(n_{\phi(1)},\dotsc, n_{\phi(m-1)},n_m)}, \dotsc$$ Since $G$ is not a winning strategy, $\{\vc_{(n_1)},\dotsc, \vc_{(n_{\phi(1)},\dotsc, n_{\phi(m-1)},n_m)},\dotsc\}$ witnesses the $\iobs$-Hurewicz property of $X$. For $B\in \bk$ there is a sequence $\{\delta_n:n\in \omega\}$ of positive real numbers such that $\{m:B^{\delta_m}\nsubseteq U \text{ for any } U\in \vc_{(n_{\phi(1)},\dotsc, n_{\phi(m-1)},n_m)}\}\in \ic_1$. Since $\ic_1\leq_\phi \ic_2$, $\phi^{-1}(\{m:B^{\delta_m}\nsubseteq U \text{ for any } U\in \vc_{(n_{\phi(1)},\dotsc, n_{\phi(m-1)},n_m)}\})\in \ic_2$. We show that $\{\vc_{n_{\phi(1)}},\dotsc, \vc_{(n_{\phi(1)},\dotsc ,n_{\phi(m)})},\dotsc\}$ witnesses the $\itbs$-Hurewicz property. To show this, for $B\in \bk$ we need to find a sequence $\{\sigma_m:m\in \omega\}$ of positive real numbers such that $\{m:B^{\sigma_m}\nsubseteq U  \text{ for any } U\in \vc_{(n_{\phi(1)},\dotsc ,n_{\phi(m)})}\}\in \ic_2$. Define $\sigma_m=\delta_{\phi(m)}$ for each $m$. Now it is easy to see that $\{m:B^{\sigma_m}\nsubseteq U  \text{ for any } U\in \vc_{(n_{\phi(1)},\dotsc ,n_{\phi(m)})}\}\subseteq \phi^{-1}(\{m:B^{\delta_m}\nsubseteq U \text{ for any } U\in \vc_{(n_{\phi(1)},\dotsc, n_{\phi(m-1)},n_m)}\}\in \ic_2$. Hence $F$ is not a winning strategy in the $\itbs$-Hurewicz game.
\ep

The class of ideals having a pseudounion can also be used to characterize the $\bs$-Hurewicz property and the corresponding game.
\begin{Th}
\label{TIH-3}
Let $\bk$ be a bornology on a metric space $X$ with a closed base. Let the ideal $\ic$ on $\omega$ have a pseudounion. The following statements are equivalent.\\
\noindent$(1)$ $X$ has the $\ibs$-Hurewicz property.\\
\noindent$(2)$ $X$ has the $\bs$-Hurewicz property.
\end{Th}
\bp
We prove only $(1)\Rightarrow (2)$. Let $A$ be a pseudounion of $\ic$. Let $C=\omega\setminus A$ and $\sigma:A\rightarrow C$ be an injection. Let $\{\uc_n:n\in \omega\}$ be a sequence of open $\bs$-covers of $X$ and write $\uc_n=\{U_{n,m}:m\in \omega\}$ for each $n$. Let $\vc_n=\{U_{1,m_1}\cap \dotsc \cap U_{n,m_n}:U_{1,m_1}\in \uc_1,\dotsc,U_{n,m_n}\in \uc_n\}$ for each $n$. Clearly $\{\vc_n:n\in \omega\}$ is a sequence of open $\bs$-covers of $X$. By $(1)$, choose a sequence $\{\wc_n:n\in \omega\}$ of finite sets $\wc_n\subseteq \vc_n$ such that for $B\in \bk$ there is a sequence $\{\delta_n:n\in \omega\}$ of positive real numbers satisfying $\{n:B^{\delta_n}\nsubseteq V \text{ for any } V\in \wc_n\} (=K)\in \ic$. Define $k_n=n$ if $n\in C$ and $k_n=\sigma(n)$ if $n\in A$. Consider the sequence $\{\wc_{k_n}:n\in \omega\}$. We show that for a $B\in \bk$ there is a sequence $\{\varepsilon_n:n\in \omega\}$ of positive real numbers such that the set $L=\{n:B^{\varepsilon_n}\nsubseteq V \text{ for any } V\in \wc_{k_n}\}$ is finite. Define $\varepsilon_n= \delta_n$ if $n\in C$ and $\varepsilon_n=\delta_{\sigma(n)}$ if $n\in A$. We have $K\subseteq^\ast A$. Consequently $K\setminus A$ and $C\cap K$ are finite and $C\cap L=C\cap K$. When $n\in A$ and $B^{\varepsilon_n}\nsubseteq V$ for any $V\in \wc_{k_n}$, we have $B^{\delta_{\sigma(n)}}\nsubseteq V$ for any $V\in \wc_{\sigma(n)}$ (as $k_n=\sigma(n)$ and $\varepsilon_n=\delta_{\sigma(n)}$). Therefore $\sigma(n)\in K$ and so $n\in A\cap L$ implies that $\sigma(n)\in C\cap K$. This shows that $A\cap L$ is finite and hence $L$ must be finite. De-constructing $\wc_{k_n}$ let $\zc_i$ be the collection of all the $i$th components of members of $\wc_{k_1}$ for $i\leq k_1$. Let $\zc_i$ be the collection of all the $i$th components of members of $\wc_{k_n}$ for $k_{n-1}<i\leq k_n$, $n>1$. Clearly $\zc_i$ is a finite subset of $\uc_i$ for each $i$. The proof will be complete if we can show that $\{\zc_i:i\in \omega\}$ witnesses the $\bs$-Hurewicz property. Let $B\in \bk$. We need to find a sequence $\{\alpha_i:i\in \omega\}$ of positive real numbers such that $\{i:B^{\alpha_i}\nsubseteq U \text{ for any } U\in \zc_i\}(=S)$ is finite. Choose $\alpha_i=\varepsilon_1$ if $i\leq k_1$ and $\alpha_i=\varepsilon_n$ if $k_{n-1}<i\leq k_n$. By construction of $\zc_i$ for $i\in \omega$ there is a $n$ satisfying $k_{n-1}<i\leq k_n$ and for $U\in \zc_i$ there is a $V\in \wc_{k_n}$ satisfying $V\subseteq U$. Now when $i\in S$ we have $B^{\alpha_i}\nsubseteq U$ for any $U\in \zc_i$ and this implies that $B^{\varepsilon_n}\nsubseteq V$ for any $V\in \wc_{k_n}$. Therefore $n\in L$. Since $L$ is finite, $S$ must be finite. Hence $\{\zc_i:i\in \omega\}$ witnesses the $\bs$-Hurewicz property.
\ep

\begin{Th}
\label{TIH-4}
Let $\bk$ be a bornology on a metric space $X$ with a closed base. Let the ideal $\ic$ on $\omega$ have a pseudounion. The following statements are equivalent.\\
\noindent$(1)$ ONE has no winning strategy in the $\ibs$-Hurewicz game.\\
\noindent$(2)$ ONE has no winning strategy in the $\bs$-Hurewicz game.
\end{Th}
\bp
We prove only $(1)\Rightarrow (2)$. Let $F$ be a strategy for ONE in the $\bs$-Hurewicz game. Define a strategy $G$ for ONE in the $\ibs$-Hurewicz game as follows. Let $F(\emptyset)=\{U_{(n)}:n\in \omega\}$, an open $\bs$-cover of $X$. Define $G(\emptyset)=\{U_{(n)}:n\in \omega\}$. Suppose that $U_{(n_1,\dotsc,n_{m-1})}$ is chosen. Let $F(\wc_{(n_1,\dotsc,n_{m-1})})=\{U_{(n_1,\dotsc, n_m)}:m\in \omega\}$. Define $G(\vc_{(n_1,\dotsc,n_m)})=\{U_{(n)}\cap \dotsc \cap U_{(n_1,\dotsc, n_{m-1},n)}:n\in \omega\}$. Let TWO choose a finite set $\vc_{(n_1,\dotsc,n_m)}$. Let TWO choose $\wc_{(n_1,\dotsc,n_m)}$, where $\wc_{(n_1,\dotsc,n_m)}$ is chosen as follows. If $m\in C$, $\wc_{(n_1,\dotsc,n_m)}$ is the collection of $m$th components of members of $\vc_{(n_1,\dotsc,n_m)}$. If $m\in A$, $\wc_{(n_1,\dotsc,n_m)}$ is the collection of $U_{(n_1,\dotsc, n_{\sigma(m)})}$, where $U_{(n_1,\dotsc, n_m)}$ is the $m$th component of members of $\vc_{(n_1,\dotsc,n_m)}$. Continuing this way we define the strategy $G$. A play in the $\bs$-Hurewicz game is $$F(\emptyset), \wc_{(n_1)}, \dotsc,F(\wc_{(n_1,\dotsc,n_{m-1})}),\wc_{(n_1,\dotsc,n_m)}, \dotsc$$ Correspondingly a play in the $\ibs$-Hurewicz game is $$G(\emptyset), \vc_{(n_1)},\dotsc, G(\vc_{(n_1,\dotsc,n_{m-1})}),\vc_{(n_1,\dotsc,n_m)}, \dotsc$$ Since $G$ is not a winning strategy, $\{\vc_{(n_1)},\dotsc, \vc_{(n_1,\dotsc,n_m)},\dotsc\}$ witnesses the $\ibs$-Hurewicz property. We show that $\{\wc_{(n_1)},\dotsc, \wc_{(n_1,\dotsc,n_m)},\dotsc\}$ witnesses the $\bs$-Hurewicz property. Let $B\in \bk$. There is a sequence $\{\delta_m:m\in \omega\}$ of positive real numbers such that $\{m:B^{\delta_m}\nsubseteq V \text{ for any } V\in \vc_{(n_1,\dotsc,n_m)}\}(=K)\in \ic$. Define $\varepsilon_m=\delta_m$ if $m\in C$ and $\varepsilon_m=\delta_{\sigma(m)}$ if $m\in A$. Consider $L=\{m:B^{\varepsilon_m}\nsubseteq W \text{ for any } U\in \wc_{(n_1,\dotsc,n_m)}\}$. We have $K\subseteq^\star A$. It can be shown that $A\cap L$ and $C\cap L$ are finite. Consequently $L$ is finite. Therefore $\{\wc_{(n_1)},\dotsc, \wc_{(n_1,\dotsc,n_m)},\dotsc\}$ witnesses the $\bs$-Hurewicz property. Hence $F$ is not a winning strategy for ONE in the $\bs$-Hurewicz game.
\ep

Combining Theorem \ref{TIH-3}, Theorem \ref{TIH-4}, \cite[Theorem 4.2]{dcpdsd}, \cite[Theorem 3.7]{sd} and \cite[Proposition 3.1]{sd} we have the following result.
\begin{Th}
\label{TIH-5}
Let $\bk$ be a bornology on a metric space $X$ with a closed base. Let the ideal $\ic$ on $\omega$ have a pseudounion. The following statements are equivalent.\\
\noindent$(1)$ $X$ has the $\ibs$-Hurewicz property.\\
\noindent$(2)$ $X$ has the $\bs$-Hurewicz property.\\
\noindent$(3)$ $X$ satisfies $\Sf(\obs,\ic\text{-}\obs^{gp})$.\\
\noindent$(4)$ $X$ satisfies $\Sf(\obs,\obs^{gp})$.\\
\noindent$(5)$ ONE has no winning strategy in the $\ibs$-Hurewicz game.\\
\noindent$(6)$ ONE has no winning strategy in the $\bs$-Hurewicz game.\\
\noindent$(7)$ ONE has no winning strategy in the game $\Gf(\obs,\ic\text{-}\obs^{gp})$.\\
\noindent$(8)$ ONE has no winning strategy in the game $\Gf(\obs,\obs^{gp})$.
\end{Th}

When the ideal $\ic$ has a pseudounion, the above results can be represented in Figure \ref{diag6}.
\begin{figure}[h]
 \begin{adjustbox}{keepaspectratio,center}
\begin{tikzcd}[column sep=2ex,row sep=5ex,arrows={crossing over}]
I\nuparrow \bs\text{-}Hurewicz\text{ }game\arrow[r,thick,equal]&I\nuparrow \ibs\text{-}Hurewicz \text{ }game\arrow[r,thick,equal]&I\nuparrow \Gf(\obs,\ic\text{-}\obs^{gp})\arrow[r,thick,equal]&I\nuparrow \Gf(\obs,\obs^{gp})&\\
\bs\text{-}Hurewicz\text{ } property\arrow[r,thick,equal]\arrow[u,thick,equal]& \ibs\text{-}Hurewicz\text{ } property\arrow[r,thick,equal]\arrow[u,thick,equal]&\Sf(\obs,\ic\text{-}\obs^{gp})\arrow[r,thick,equal]\arrow[u,thick,equal]&\Sf(\obs,\obs^{gp})\arrow[u,thick,equal]&\\
\end{tikzcd}
\end{adjustbox}
 \caption{Diagram of the selection principles when $\ic$ has a pseudounion}
 \label{diag6}
\end{figure}

{}

\end{document}